\documentclass[11pt,letterpaper]{article}
\usepackage{ifthen,latexsym,amssymb,amsmath}

\newcommand{\B}[1]{{\bf #1}}
\newcommand{\I}[1]{{\mathbb #1}}
\newcommand{\OO}[1]{\overline{#1}}

\newcommand{\e}{\varepsilon}

\renewcommand{\mid}{:}

\newcommand{\beq}[1]{\begin{equation}\label{eq:#1}}
\newcommand{\eeq}{\end{equation}}
\newcommand{\req}[1]{\textrm{(\ref{eq:#1})}}

\newtheorem{theorem}{Theorem}
\newcommand{\bth}[2][nothing]{\ifthenelse{\equal{#1}{nothing}}
 {\begin{theorem}} {\begin{theorem}[#1]}\label{th:#2}}

\newtheorem{lemma}[theorem]{Lemma}
\newcommand{\blm}[2][nothing]{\ifthenelse{\equal{#1}{nothing}}
 {\begin{lemma}} {\begin{lemma}[#1]}\label{lm:#2}}

\newtheorem{problem}[theorem]{Problem}
\newcommand{\bpr}[2][nothing]{\ifthenelse{\equal{#1}{nothing}}
 {\begin{problem}} {\begin{problem}[#1]}\label{pr:#2}}

\newtheorem{conjecture}[theorem]{Conjecture}

\newtheorem{claim}[theorem]{Claim}

\newcommand{\Claim}[2]{\begin{claim}\label{cl:#1} #2\end{claim}}

\newcommand{\bpf}[1][Proof.]{\smallskip\noindent{\it #1} }
\newcommand{\qed}{\nolinebreak\mbox{\hspace{5 true pt}%
  \rule[-0.85 true pt]{3.9 true pt}{8.1 true pt}}}
\newcommand{\cqed}{\nolinebreak\mbox{\hspace{5 true pt}%
  \rule[-0.85 true pt]{2.0 true pt}{8.1 true pt}}}
\newcommand{\epf}{\qed \medskip}
\newcommand{\bcpf}{\bpf[Proof of Claim.]}
\newcommand{\ecpf}{\cqed \medskip}

\newcommand{\brm}{\smallskip\noindent{\bf Remark.} }

\newcommand{\V}[1]{{\bf #1}}
\newcommand{\blow}[2]{#1(\!(#2)\!)}
\newcommand{\grad}{\mathrm{grad}}
\newcommand{\rn}{\boldsymbol}
\newcommand{\pp}{P}
\newcommand{\barL}{\OO L}

\title{Minimum Number of $k$-Cliques in Graphs with Bounded Independence Number}
\author{Oleg Pikhurko\footnote{Supported by  the
European Research Council (grant agreement no.~306493)
and National Science Foundation of the USA (grant DMS-1100215).}\\
Mathematics Institute and DIMAP\\
University of Warwick\\
Coventry CV4 7AL, UK
\and
Emil R.\ Vaughan\\
Centre for Discrete
Mathematics\\ 
Queen Mary University of London\\
London E1 4NS, UK
}
\date{}

\begin{document}

\maketitle

\begin{abstract}
Erd\H os asked in 1962 about the value of $f(n,k,l)$, the minimum
number of $k$-cliques in a graph with order $n$ and independence number 
less than $l$. The case $(k,l)=(3,3)$ was solved by Lorden.
Here we solve the problem (for all large $n$) for $(3,l)$ with $4\le l\le 7$ and $(k,3)$ with $4\le k\le 7$. Independently,
Das, Huang, Ma, Naves, and Sudakov
resolved the cases $(k,l)=(3,4)$ and $(4,3)$. 
\end{abstract}

\section{Introduction}

Let us give some definitions first. As usual, a \emph{graph} $G$ is a pair
$(V(G),E(G))$, where $V(G)$ is the \emph{vertex set} and the \emph{edge
set} $E(G)$ consists of unordered pairs of vertices. An \emph{isomorphism}
between graphs $G$ and $H$ is a bijection $f:V(G)\to V(H)$ that preserves
edges and non-edges.  For a graph $G$, let $\OO G=(V(G),{V(G)\choose
2}\setminus E(G))$
denote its \emph{complement} and $v(G)=|V(G)|$
denote its \emph{order}. For graphs $F$ and $G$
with 
$v(F)\le v(G)$, let $\pp(F,G)$ be the number of $v(F)$-subsets of $V(G)$
that induce in $G$ a subgraph isomorphic to $F$; further, define the 
\emph{density of $F$ in $G$} to be
 \beq{density}
 p(F,G)=\pp(F,G)\,{v(G)\choose v(F)}^{-1}.
 \eeq

Let $K_k$ denote the
complete graph on $k$
vertices. Let $\alpha(G)=\max\{l\mid \pp(\OO K_l,G)>0\}$
be the \emph{independence number} of $G$, that is, the maximum size
of an edge-free set of vertices.

Given a graph $F$ on
$[m]=\{1,\dots,m\}$ and a sequence of disjoint
sets $V_1,\dots,V_m$,
let the \emph{expansion} 
$\blow{F}{V_1,\dots,V_m}$ be the graph on $V_1\cup\dots\cup V_m$ obtained
by
putting the complete graph on each $V_i$ and
putting, for each edge $\{i,j\}\in E(F)$, the complete bipartite graph between
$V_i$ and $V_j$. An expansion is
\emph{uniform} if  $\big|\,
|V_i|-|V_j|\,\big|\le 1$
for any $i,j\in [m]$.  
If we consider expansion in terms of complements, then it amounts
to blowing up each vertex $i$ of $\OO F$ by factor $n_i$ (and taking the
complement of the obtained graph). Clearly, expansions cannot increase the
independence
number.

We consider the following extremal function
 $$
 f(n,k,l)=\min\left\{\pp(K_k,G)\mid v(G)=n,\, \alpha(G)<l\right\},
 $$
 that is, the minimum number of $k$-cliques in a graph with $n$
vertices that does not contain $\OO K_l$. 
This function (in its full generality) was first defined 
by Erd\H os~\cite{erdos:62a} in 1962.

Earlier, Goodman~\cite{goodman:59} determined $f(2n,3,3)$; his bounds
also give the asymptotic value of $f(2n+1,3,3)$. 
 Lorden~\cite{lorden:62}
determined $f(n,3,3)$ and showed that the complement of
$T_2(n)$ is the unique extremal graph when $n\ge 12$, where
the \emph{Tur\'an graph} $T_m(n)$ is the complete $m$-partite
graph on $[n]$ with parts being nearly equal. (In other words,
$T_m(n)$ is the complement of the
uniform expansion of $\OO K_m$.)

Erd\H os~\cite{erdos:62a} asked if perhaps 
 \beq{E62a}
 f(n,k,l)=\pp(K_k,\OO T_{l-1}(n)),
 \eeq
 that is, if the uniform expansion of $\OO K_{l-1}$ gives the value of
$f(n,k,l)$ and,
specifically, if 
 \beq{E62b}
 f(3n,3,4)=3 {n\choose 3}.
 \eeq
 Nikiforov~\cite{nikiforov:01} showed that the limit
 \beq{ckl}
 c_{k,l}=\lim_{n\to\infty} \frac{f(n,k,l)}{{n\choose k}}
 \eeq
 exists for every pair $(k,l)$ and  that the lower
bound $c_{k,l}\ge (l-1)^{1-k}$ given by the graphs 
$\OO T_{l-1}(n)$ as $n\to\infty$ can be sharp only for finitely
many pairs $(k,l)$. Thus, it was too optimistic to expect
that \req{E62a} holds.

The main motivation of the papers \cite{erdos:62a,goodman:59} came from
Ramsey's theorem~\cite{ramsey:30}, which implies that $f(n,k,l)>0$
when $n\ge n_0(k,l)$ is sufficiently large. Both papers also
considered the related problem of minimising $p(K_k,G)+p(\OO K_k,G)$
over an (arbitrary) order-$n$ graph $G$. The last question, known as the
\emph{Ramsey multiplicity problem}, attracted a lot of attention
and led to many important developments.

On the other hand, the problem of 
determining $f(n,k,l)$ was rather neglected although it was
mentioned in Bollob\'as' book \cite[Problem 11 on Page
361]{bollobas:egt} and Thomason's survey~\cite[Section
5.5]{thomason:02}. One possible reason is that determining
$c_{k,l}$, even for some small $k$ and $l$, 
might require keeping track of too many different subgraph densities
than what is practically feasible when doing calculations ``by hand''.

Razborov~\cite{razborov:07} introduced a powerful 
formal system for deriving inequalities
between subgraph densities, where a computer can be
used to do routine book-keeping. One aspect of his theory
(introduced in \cite{razborov:10})
allows us to minimise linear combinations 
of subgraph densities by setting up and solving
a semi-definite program. In some cases, the
numerical solution thus obtained can be converted into a rigorous
mathematical proof. Baber and Talbot~\cite{baber+talbot:11}
and Vaughan~\cite{flagmatic:2.0} 
(see~\cite{falgas+vaughan:12,falgas+vaughan:13})
wrote openly available software for doing such calculations.

By using \emph{Flagmatic}~\cite{flagmatic:2.0}, we can solve
the problem (for all large $n$) when $k=3$ with $4\le l\le 7$ or $l=3$ with
$4\le k\le 7$.
Independently, Das, Huang, Ma, 
Naves, and Sudakov~\cite{das+huang+ma+naves+sudakov:12:arxiv} solved the problem
when $n$ is large and $(k,l)=(3,4)$ or $(4,3)$, also by
using flag algebras.

We state our results as three separate theorems.

\begin{theorem}[Asymptotic Result]\label{th:asympt}
 \begin{eqnarray}
 c_{3,l} &=& (l-1)^{-2},\qquad 4\le l\le 7,\label{eq:c34}\\
 c_{4,3} &=& 3/25,\label{eq:c43}\\
 c_{5,3} &=& 31/5^4\ =\ 31/625,\label{eq:c53}\\
 c_{6,3} & =& 19211/2^{20}\ =\ 19211/1048576,\label{eq:c63}\\
 c_{7,3} &=& 98491/2^{24}\ =\ 98491/16777216.\label{eq:c73}
 \end{eqnarray}
 Furthermore, we have in each of these cases that
 \beq{asympt}
 f(n,k,l)=c_{k,l}{n\choose k}+O(n^{k-1}).
 \eeq
 \end{theorem}

The upper bounds in \req{c34}, \req{c43}, and \req{c53} are obtained
by taking a uniform expansion of $F$, where $F$ is respectively $\OO
K_{l-1}$, the $5$-cycle
$C_5$, and (again) $C_5$. Easy calculations show that the
density of $k$-cliques in these graphs is as required.
These upper bounds on $c_{4,3}$ and $c_{5,3}$ come from 
Nikiforov~\cite{nikiforov:01}.
In a subsequent paper~\cite{nikiforov:05:arxiv}, he also showed
that an order-$n$ graph $G$ with $\alpha(G)<3$ satisfies
$P(K_4,G)\ge (\frac{3}{25}+o(1)) {n\choose 4}$ under the
additional assumption that $G$ is close to being regular.

The upper bounds in \req{c63} and \req{c73} come from a more
complicated construction. The \emph{Clebsch
graph} 
$L$ has binary $5$-sequences of even weight (i.e.\ with an even number
of entries equal to 1) for  vertices, with
two vertices being adjacent
if the term-wise sum modulo 2 of the corresponding sequences has weight $4$. 
For example, the neighbours of $00011\in V(L)$ are $01100$, $10100$, $11000$, $11101$, and $11110$.
It easily follows from this description that the
Clebsch graph is triangle-free and vertex-transitive.
For example, an automorphism that maps $00000$ to $11000$ is to flip the
first two bits.

The complement $F=\barL$ of the Clebsch graph is a
10-regular graph on $16$ vertices.  Take a uniform expansion $F'$ of $F$
of large order $n$.
The limit of $p(K_k,F')$ as $n\to\infty$ is equal to the probability that, if we sample independent
uniformly distributed vertices $x_1,\dots,x_k\in V(L)$, they do not induce any edge
in $L$. 
By the vertex-transitivity of $L$, we can fix $x_1=00000$. The Clebsch graph
has the following maximal independent sets containing $00000$: the 
sequences that we add to $00000$ must have weight 2, with the corresponding
pairs of indices forming either $K_{1,4}$ (the star with 4 edges)
or $K_3$ (the triangle). There are 5 of the former sets and $10$ of the latter
sets, of sizes 5 and 4 respectively. A
straightforward inclusion-exclusion counting shows that the above probability is
 $$
 \frac{5\cdot 5^{k-1} + 10\cdot 4^{k-1} - 30\cdot 3^{k-1} + 20\cdot 2^{k-1}-4}{16^{k-1}}.
 $$
 By plugging in $k=6$ and $7$, we get the upper bounds on $c_{k,3}$ stated in
\req{c63} and
\req{c73}. 

The upper bound in~\req{asympt} follows by observing that if
we pick a random injection $\rn{\phi}:[k]\to V(F')$, where $F'$ a uniform
expansion of $F$ of order $n$, and condition on the restriction of $\rn{\phi}$ 
to $[i]$ for $i<k$, then the probability
that $\rn{\phi}(i+1)$ belongs to 
a particular part of $F'$ is $1/v(F)+O(1/n)$. Thus $p(K_k,F')$ is within
additive term $O(1/n)$ from its limit as $n\to\infty$.

The lower bounds of Theorem~\ref{th:asympt} are proved in Section~\ref{lower} by using flag algebras.

We say that two graphs $G$ and $H$ of the same order are at 
\emph{edit distance at most $m$} or \emph{$m$-close} if $G$
can be made isomorphic to $H$ by changing (adding or deleting) at most $m$
edges. By inspecting the proof certificate
returned by a flag algebra
computation, one can sometimes describe
the structure of all almost extremal graphs up to a small edit distance
(see, for
example, \cite{CKPSTY,hatami+hladky+kral+norine+razborov:13,pikhurko:11}). This
also works
here and we can establish the following results that
apply when $(k,l)$ is one of the pairs $(3,l)$ with $3\le l\le 7$, $(k,3)$ with $4\le k\le 5$,
and $(k,3)$ with $6\le k\le 7$, while $F$ is respectively $\OO K_{l-1}$, $C_5$, and
$\barL$.

\begin{theorem}[Stability Property]\label{th:stab}  Let $k,l,F$ be as
above. Then for
every 
$\e>0$ there exist $\delta>0$ and $n_0$ such that every
graph $G$ of order $n\ge n_0$ with $\alpha(G)<l$ and
$P(K_k,G)\le (c_{k,l}+\delta){n\choose k}$ is $\e {n\choose 2}$-close
to a uniform expansion of $F$.\end{theorem}

We see that, in each case above, almost extremal graphs on $[n]$ have 
the same structure  up to the edit
distance of $o(n^2)$. Such extremal problems are called \emph{stable}. The stability property, besides being
of interest on its own, is often very helpful in establishing the
exact result for all large $n$. Here, we also use stability to prove the
following theorem.

\begin{theorem}[Exact Result]\label{th:exact}   Let $k,l,F$ be
as above. Then there is
$n_0$ such that every graph $G$ of order $n\ge n_0$ with $\alpha(G)<l$ and
the minimum number of $K_k$-subgraphs contains an expansion
$F'=\blow{F}{V_1,\dots,V_m}$ as a spanning subgraph 
(that is, $V_1\cup\dots\cup V_m=V(G)$ and $E(F')\subseteq E(G)$). 

\end{theorem}

Let $n$ be sufficiently large. Since $G$ in Theorem~\ref{th:exact} is extremal and $F'$ is $\OO K_{l-1}$-free, we have that
$\pp(K_k,G)=\pp(K_k,F')$, that is, the value of $f(n,k,l)$ 
is attained
by some expansion of $F$. Furthermore, if $l=3$ and $4\le k\le 7$, then
$G$ is necessarily
equal to $F'$  because the addition of any
extra edge to $F'$ creates at least one copy of $K_k$. Next, consider the four remaining
cases, that is, $k=3$ and $4\le l\le 7$. It is easy to show that $\OO T_{l-1}(n)$ has the 
smallest number of triangles among all
order-$n$ expansions of   $\OO K_{l-1}$. Thus 
Theorem~\ref{th:exact} proves Erd\H os' conjecture~\req{E62b} for
all large $n$. However note that there are other extremal constructions for  $f(n,3,l)$ with $4\le l\le 7$ that can be
obtained from $\OO T_{l-1}(n)$ by adding edges so that no new triangles are
created.

It would
be interesting to determine those $l$ for which $c_{3,l}=(l-1)^{-2}$. We know now that this is the case for all $2\le l\le 7$. Nikiforov~\cite{nikiforov:01} showed that this equality can hold for only finitely many $l$. 
Das et al~\cite{das+huang+ma+naves+sudakov:12:arxiv} proved that no $l\ge 2074$
satisfies it.

Although our proofs rely on extensive computer calculations, new
mathematical ideas are also introduced (such as, for example,
Theorem~\ref{th:1} that deals with all studied cases in a unified manner).
Hopefully, these ideas and results will be useful for other problems.
For example, the concept of a \emph{phantom edge} introduced here 
in Section~\ref{phantom} has been successfully applied to another extremal problem~\cite{falgas+marchant+pikhurko+vaughan}.

\section{Notation}

Here we collect some graph theory notation that we use. 

The cycle (resp.\ path) with $k$ vertices is denoted by $C_k$ (resp.\ $P_k$).

Let $G$ and $H$ be graphs. We write $H\subseteq G$ and say that $H$ is a \emph{subgraph} of $G$ if $V(H)\subseteq V(G)$ and $E(H)\subseteq E(G)$. A
subgraph $H\subseteq G$ is called \emph{spanning} if $V(H)=V(G)$. It is called 
\emph{induced}  if $H=G[\,V(H)\,]$, where we denote $G[X]=(X,\{\{x,y\}\in E(G)\mid x,y\in X\})$ for $X\subseteq V(G)$. A \emph{strong homomorphism} from $H$ to $G$ is a map $\phi: V(H)\to
V(G)$ that preserves both edges and non-edges. 
For example, $H$ admits a strong homomorphism to $K_2$ if and only if $H$ is a complete bipartite graph. An \emph{embedding} is a strong homomorphism which is injective; in other words, it is an isomorphism from $H$ to an induced subgraph of $G$.

An \emph{automorphism} of $G$ is a map $V(G)\to V(G)$ that preserves both edges and non-edges (i.e.\  an isomorphism
of $G$ to itself).
A graph $G$ is \emph{vertex-transitive} if for every two vertices
there is an automorphism of $G$ mapping one into the other.
The \emph{neighbourhood} of a vertex $x\in V(G)$ is 
 $$
 \Gamma_G(x)=\big\{y\in V(G)\mid \{x,y\}\in E(G)\big\}.
 $$
 The \emph{closed neighbourhood} of $x$ is $\hat
\Gamma_G(x)=\Gamma_G(x)\cup\{x\}$.

The \emph{Ramsey number} $R(k,l)$ is the minimum $n$ such that every order-$n$ graph
has a $k$-clique or an independent set of size $l$. Thus $f(n,k,l)>0$ if and only if $R(k,l)\ge n$.

\section{Lower Bounds in Theorem~\ref{th:asympt}}\label{lower}

\subsection{Proof Certificates}

As we have already mentioned, our lower bounds are proved with the help of
a computer 
by using flag algebras and semi-definite programming, see
Razborov~\cite{razborov:07,razborov:10}. This method 
is described in a number of research publications
(\cite{baber+talbot:11,falgas+vaughan:12,
falgas+vaughan:13,hirst:4vertex,razborov:07,razborov:10}), so we will be brief.

We used \emph{Flagmatic} (Version 2.0) \cite{flagmatic:2.0} for the computations.
For each proof we present, we provide a certificate that contains the information
needed for others to be able to verify all claims. The script \texttt{inspect\_certificate.py}
that comes with \emph{Flagmatic} can be used for investigating the certificates and performing
some level of verification. The certificates are in a documented format \cite{flagmatic:2.0} and it
is hoped that others will be able to independently verify them.

Also, we include the code that generated each certificate as well as
the transcript of each session, 
to aid the reader in repeating our calculations. This may be helpful if the
reader would like to experiment with the software by changing parameters (or to
apply \emph{Flagmatic} to some related problems).

These materials are available from 
\emph{Flagmatic}'s website at
 $$
 \mbox{\texttt{http://flagmatic.org/examples/Fkl.tgz}}
 $$
 Each solved case $(k,l)$ is supported by the following data: the complete flagmatic code,
the transcript of the session, and all generated certificates. For example, the corresponding files for the case $(k,l)=(7,3)$ are \texttt{73.sage}, \texttt{73.txt}, and two certificates \texttt{73.js} and \texttt{73a.js}.

Alternatively, the ancillary folder of~\cite{pikhurko+vaughan:Fkl:arxiv}
contains all files
except some certificates whose sizes are larger than arxiv's
allowance. The reader should be able to generate these certificates by running the appropriate scripts
with \emph{Flagmatic
2.0}.

Also, the cases $(3,4)$ and $(k,3)$ with $4\le k\le 7$ were previously solved
with Version 1.5 of \emph{Flagmatic};
see~\cite{pikhurko+vaughan:Fkl:arxiv} (Version 3) for all details. This is
reassuring as \emph{Flagmatic 2.0} was re-written essentially from scratch (when it was
decided to do everything inside \emph{sage} for greater functionality).

Our presentation is different from that of Das et
al~\cite{das+huang+ma+naves+sudakov:12:arxiv} who worked hard on
making their paper self-contained and the proof
as human-readable as possible. This has many advantages (such as giving more insight into
the problem) but makes the paper rather long. Our objective is to
present
formal rigorous proofs of all claimed results. We do so by describing the
information that is contained in the certificates and by showing how it implies
the stated results. While the certificates are not very suitable for direct
inspection (some of them are very large and contain integers with hundreds 
of digits), the reader may verify all stated properties by using
\emph{Flagmatic} or by writing an independent script. 

Let us give some definitions that are needed to describe the certificates. Fix one
of the pairs $(k,l)$ as above.

Let us call a graph \emph{admissible} if its independence number is less than
$l$.
A \emph{type} is a pair $(H,\phi)$ where $H$ is an admissible
graph and $\phi:[v]\to V(H)$ is a bijection, where $v=v(H)$. 
Given a type
$\tau=(H,\phi)$ as above,
a \emph{$\tau$-flag} is a pair $(G,\psi)$ where $G$ is an admissible
graph and $\psi:[v]\to V(G)$
is an injection such that $\psi\circ \phi^{-1}: V(H)\to V(G)$ 
is an embedding (that is, an injection that
preserves both edges and non-edges). Informally, a type is a
vertex-labelled graph and a
$\tau$-flag is
a partially labelled graph such that the labelled vertices induce 
$\tau$. The \emph{order} $v((G,\psi))$ of a type or a flag 
is $v(G)$, the number of vertices in it.

For two $\tau$-flags $(G_1,\psi_1)$ and $(G_2,\psi_2)$
with $n_1\le n_2$ vertices, let $\pp((G_1,\psi_1),(G_2,\psi_2))$
be the number of $n_1$-subsets $X\subseteq V(G_2)$ such that
$X\supseteq \psi_2([v])$ (i.e.\ $X$ contains all labelled 
vertices) and the $\tau$-flags 
$(G_1,\psi_1)$ and $(G_2[X],\psi_2)$ are \emph{isomorphic},
meaning that there is a graph isomorphism that preserves
the labels. Also, define the \emph{density}
 $$
 p((G_1,\psi_1),(G_2,\psi_2))=\frac{\pp((G_1,\psi_1),(G_2,\psi_2))}{{
n_2-v\choose n_1-v}},
 $$
 to be the probability that a uniformly drawn random $n_1$-subset $X$ of
$V(G_2)$ with 
$X\supseteq \phi_2([v])$ induces a copy of the $\tau$-flag $(G_1,\psi_1)$ in
$(G_2,\psi_2)$.

Now, we can present the information that is contained
in each certificate (a file with extension \texttt{js}) and is needed in the
proof.

First, the certificate lists all (up to an isomorphism) admissible $N$-vertex graphs for
some integer 
$N$.
Let us denote these graphs by $G_1,\dots,G_g$.
Then the certificate describes some types $\tau_1,\dots,\tau_t$ such that their
graph components are pairwise
non-isomorphic (as unlabelled graphs) and 
$N-v(\tau_i)$ is a positive even number for each $i\in [t]$.

The certificate contains, for each $i\in [t]$,  the list
$(F_1^{\tau_i},\dots,F_{g_i}^{\tau_i})$ of all $\tau_i$-flags (up to
isomorphism of $\tau_i$-flags) with exactly
$(N+v(\tau_i))/2$ vertices. 

Also, for each $i\in [t]$, the certificate (indirectly) contains a symmetric
positive semi-definite $g_i\times g_i$-matrix $Q^{\tau_i}$. More precisely, 
the matrix $Q^{\tau_i}$ is represented in the following manner: we have a
diagonal matrix $Q'$ all whose diagonal entries are positive rational
numbers and a rational matrix $R$ such that 
 \beq{R}
 Q^{\tau_i}=RQ'R^T.
 \eeq This
decomposition automatically implies that the matrix $Q^{\tau_i}$ is
positive
semi-definite.

Now, let $G$ be an admissible graph of large order $n$. Initially, let $a=0$.
Let
us do the
following for each $v$ such that $N-v$ is a positive even integer. Enumerate all
$n(n-1)\dots (n-v+1)$ injections $\psi:[v]\to V(G)$. If the \emph{induced type}  $G[\psi]=(G[\,\psi([v])\,],\psi)$ is isomorphic to some $\tau_i$ (as
vertex-labelled graphs), then we add $\V x_\psi
Q^{\tau_i}\V x_\psi^T$ to
$a$, where
 \beq{x}
 \V
x_\psi=\big(\pp(F_1^{\tau_i},(G,\psi)),\dots,\pp(F_{g_i}^{\tau_i},(G,
\psi))\big).
 \eeq
 Since each $Q^{\tau_i}$ is positive semi-definite, we have that $\V x_\psi
Q^{\tau_i}\V x_\psi^T\ge 0$ and that the final $a$ is non-negative. 

Let us take some type $\tau$  of order $v$ and two $\tau$-flags $F_1$ and $F_2$
with respectively $\ell_1$ and $\ell_2$ vertices. Let $\ell=\ell_1+\ell_2-v$.
Consider the sum
 \begin{equation}\label{eq:prod}
 \sum_{\psi\,:\,G[\psi]\cong \tau} \pp(F_1,(G,\psi))\, \pp(F_2,(G,\psi)),
 \end{equation} 
 taken over
all injections $\psi:[v]\to V(G)$ such that the induced type $G[\psi]$ 
is isomorphic to $\tau$. Each term
$\pp(F_i,(G,\psi))$ in (\ref{eq:prod}) can be
expanded as the sum over $\ell_i$-sets $X_i$ with
$\psi([v])\subseteq X_i \subseteq V(G)$ of the indicator
function that $(G[X_i],\psi)$ is a $\tau$-flag isomorphic to $F_i$. Ignoring the
choices
when $X_1$ and $X_2$ intersect outside of $\psi([v])$, the remaining
terms can be generated by choosing an $\ell$-set $X=X_1\cup X_2$ first, then
injective $\psi:[v]\to X$, and finally $X_1$ and $X_2$. Clearly, the terms that
we ignore contribute at most $O(n^{\ell-1})$ in total. Also, the contribution
of each $\ell$-set $X$ to \req{prod} depends only on the isomorphism class $H$ of
$G[X]$. Thus  the sum
in~(\ref{eq:prod})
can be written (modulo
an additive error term $O(n^{\ell-1})$) as an explicit
linear combination of the subgraph counts $\pp(H,G)$, where $H$ runs
over unlabelled graphs with $\ell$ vertices, see 
e.g.~\cite[Lemma~2.3]{razborov:07}. 

By the above discussion, if we expand
each quadratic form $\V x_\psi Q^{\tau_i}\V x_\psi^T$ in the definition of $a$ 
and take
the 
sum over all injections $\psi$, then we will get a
representation
 \beq{a}
 0\le a=\sum_{i=1}^g \alpha_i P(G_i,G)+O(n^{N-1}),
 \eeq
 where each $\alpha_i$ is a rational number that does not depend on $n$
and can be computed given the above
information (types, flags, and matrices). An explicit formula for
$\alpha_i$ is rather messy, so we do not state it.

The crucial property that our certificates possess is that
 \beq{main}
 \alpha_i\le p(K_k,G_i)-c_{k,l}',\quad \mbox{for every $i\in[g]$},
 \eeq
 where $c_{k,l}'$ is the right-hand side of the appropriate 
statement~\req{c34}--\req{c73}, i.e.\ $c_{k,l}'$ is the lower
bound on $c_{k,l}$ that we want to prove.
This property (involving rational numbers) can be verified by the 
stand-alone script
\texttt{inspect\_certificate.py}
that uses exact arithmetic.

If we assume that~\req{main} holds,
then we have, by Bayes' formula, that
 \beq{derivation}
 p(K_k,G)-c_{k,l}'=\sum_{i=1}^g (p(K_k,G_i)-c_{k,l}')p(G_i,G)\ge \sum_{i=1}^g
\alpha_i p(G_i,G)\ge -O(1/n).
 \eeq
 Thus we derived not only $c_{k,l}\ge c_{k,l}'$ but also the claimed lower
bound in~\req{asympt}.

At this point, we may stop and assume that Theorem~\ref{th:asympt} has been
proved (modulo verifying all the claims above with the help of a computer).
However, it may be 
useful to say a few words how these certificates were obtained. Finding matrices
$Q^{\tau_1},\dots,Q^{\tau_t}$ amounts
to solving a semi-definite program. The program is usually is quite large.
So it is generated by computer as well; \emph{Flagmatic} provides
a highly customisable way of doing this. Then the obtained program 
is fed into an SDP-solver which return floating-point matrices. 
It is a good idea to start with as small as possible $N$ and keep increasing
it until the obtained (floating-point) bound seems to be equal to the
conjectured value. We found it beneficial, at this stage, 
to use the double-precision \texttt{spda\_dd} solver that usually
returns the correct values of around 20 first decimal digits.

In fact, this was how the
extremal configuration for $c_{6,3}$ was discovered. The solver seemed to give
the same bound $c_{6,3}\ge 19211/2^{20}$ for both $N=7$ and $8$.
Here, the denominator is a high
power of 2. This suggested that an extremal configuration might be a uniform
expansion of a graph with 16 vertices which made us
to look at such graphs.

This process of converting
the obtained floating-point matrices into those
that satisfy~\req{main} exactly also uses a computer. It is fairly automated 
in \emph{Flagmatic}, although it
sometimes requires adjusting various  parameters and options. Of course,
once we have found suitable rational matrices that provide a rigorous proof, 
we can ignore their
floating-point lineage altogether.

One strategy to simplify the proof certificates once $N$ has been fixed, is to
reduce the number of types as much
as possible by re-running the SDP-solver and checking that we still get 
the same bound. Note that $\tau_1,\dots,\tau_t$ need not enumerate all types.
The removal of some type $\tau$ effectively means that we make the
corresponding matrix  $Q^{\tau}$
to be identically 0. (Likewise, $F_1^{\tau_i},\dots,F_{g_i}^{\tau_i}$ need not
enumerate all $\tau_i$-flags but this observation does not seem to be very useful.)

Another useful trick comes from the following lemma.

\begin{lemma}\label{lm:forced} Suppose that we have a flag algebra proof,
as specified above, that the value of $c_{k,l}$ is given by uniform expansions of 
a $\OO K_l$-free graph $F$. Fix $i\in [t]$. Let the $i$-th type
$\tau_i$ be $(H,\phi)$
and let $v=v(\tau_i)$.
Let $n$ be large and $G$ be a
uniform expansion of $F$ of order $n$. Let 
$\psi:[v]\to V(G)$ be an injection such that $\psi\circ \phi^{-1}$ is
an embedding of $H$ into $G$. Then
$\V x_\psi Q^{\tau_i}\V x_\psi^T=O(n^{N-v-1})$, where
$\V x_\psi$ is defined by~\req{x}.\end{lemma}

\bpf Since each part $V_i$ of $G$ is homogeneous,
any modification of the injection $\psi$ such that its values stay in the same
parts is an embedding. These new injections give the same
vector $\V x_\psi$. Thus, with $m=v(F)$,
 \beq{xQx}
 0\le \left(\frac nm+O(1)\right)^v \, \V x_\psi Q^{\tau_i}\V x_\psi^T\le a.
 \eeq

Let us run our flag algebra
proof on $G$. It shows in fact that 
$p(K_k,G)\ge c_{k,l}+a/{n\choose N}+O(1/n)$. Also, as we have
previously remarked, $p(K_k,G)$ deviates from $c_{k,l}$ by
at most $O(1/n)$. By~\req{xQx} we conclude that $a=O(n^{N-1})$, implying the lemma.\qed

Thus, when we let $n\to\infty$, the normalised limit of $\V x_\psi$ is a
zero eigenvector of $Q^{\tau_i}$. (Note that $\B xQ\B x^T=0$ for $Q\succeq 0$ implies
that $Q\B x^T=\B 0$.) We call such a zero eigenvector 
\emph{forced}. By inspecting the graph $F$
that gives the upper bound in Theorem~\ref{th:asympt}, we can identify
forced zero eigenvectors. It is crucial to know all forced zero eigenvectors during the rounding step
because a small but uncontrolled perturbation of  $Q^{\tau_i}$ may result in negative eigenvalues. \emph{Flagmatic 2.0} takes care of this by ensuring that
the column space of the matrix $R$ in~\req{R} is orthogonal to all forced zero eigenvectors of $Q^{\tau_i}$ (when an extremal construction is supplied using the function \verb=set_extremal_construction=).

Lemma~\ref{lm:forced} can be generalised to many other problems. This idea
was first used by
Razborov \cite{razborov:10}.

There are further relations that have to hold in a flag algebra proof. For $i\in [g]$, call the graph $G_i$ \emph{sharp} if
\req{main} is equality, that is, $\alpha_i=p(K_k,G_i)-c_{k,l}$. (We know by now that $c_{k,l}=c_{k,l}'$.) 

\begin{lemma}\label{lm:ForcedSharp} Suppose that we have a flag algebra proof,
as specified above, that the value of $c_{k,l}$ is given by uniform expansions of 
a $\OO K_l$-free graph $F$. Let $n$ be large and $G$ be a
uniform expansion of $F$ of order $n$. Let $i\in [g]$ be such that $G_i$ embeds into $G$. Then $G_i$ is sharp.\end{lemma}

\bpf Let $m=v(F)$. Note that $P(G_i,G)\ge (n/m+o(1))^N$: if we take an embedding $f$ of $G_i$
into $\blow{F}{U_1,\dots,U_m}$, then any injection $f':V(G_i)\to V(G)$ with $f(x)$ and $f'(x)$
belonging to the same part $U_j$ is also an embedding. 
We have by \req{a} and \req{main} that
  \begin{eqnarray}
  p(K_k,G)-c_{k,l} &\ge & 
\sum_{j\in [g]}  (p(K_k,G_j)-c_{k,l}-\alpha_j)p(G_j,G)+ O(1/n)\nonumber\\
 &\ge& (p(K_k,G_i)-c_{k,l}-\alpha_i)p(G_i,G)+ O(1/n)\label{eq:XX}\\
 &\ge& \frac{N!\,(p(K_k,G_i)-c_{k,l}-\alpha_i)}{m^N}+o(1).\nonumber
 \end{eqnarray}
  Since $p(K_k,G)-c_{k,l}=o(1)$ by our assumption, we conclude (by 
using~\req{main} again) that $G_i$ is sharp, as required.\epf

\emph{Flagmatic} also uses the restrictions given by Lemma~\ref{lm:ForcedSharp} for rounding (if a construction is provided).
In some cases, the large amount of data and/or the presence of tiny but non-zero coefficients required from us to reduce the number of types as much as
possible (essentially by trial and error) and to use the double-precision SDP-solver \verb=sdpa_dd=. 
Below we mention briefly how this process went in each solved 
case and what further actions (if any) were needed.

\subsection{Cases $(k,l)=(4,3)$ or $(5,3)$}

The rounding procedure worked without any
issues for these two cases. In both cases, we used the 6-vertex universe
that contains 38 graphs with independence number at most $2$.

\subsection{Cases $(k,l)=(6,3)$ or $(7,3)$}

In these cases, we found it more convenient to work with the complements: namely, we forbid
$K_3$ and minimise the density of $\OO K_k$ for $k=6,7$.
These cases went through without any problems. While $c_{6,3}$ could be computed
by using graphs with at most 7 vertices, it seems that
the determination of $c_{7,3}$ by this method requires 8-vertex graphs.

\subsection{Cases $k=3$ and $4\le l\le 7$}\label{phantom}

One difficulty that we had to overcome is that there are some further relations
that a flag algebra proof of $c_{3,l}\ge (l-1)^{-2}$ has to satisfy, in addition to those
given by Lemmas~\ref{lm:forced} and~\ref{lm:ForcedSharp}.

\begin{lemma}\label{lm:phantom} Suppose that we have a flag algebras proof that
$c_{3,l}\ge (l-1)^{-2}$ as above. Let $n$ be large and $T=\OO T_{l-1}(n)=\blow{\OO K_{l-1}}{V_1,\dots,V_{l-1}}$. Let $T'$ be obtained from $T$ by adding one extra edge $\{x_1,x_2\}$
between $V_1$ and $V_2$. 
If some $G_i$ admits an embedding $f$ into $T'$, then it is sharp.
\end{lemma}
 \bpf 
 Let
$\e>0$ be a small constant and let $n\to\infty$. Let the graph $G$ be obtained
from $T$ by adding all edges
between $U_1$ and $U_2$, where $U_i\subseteq V_i$ is a set of size $\lfloor \e
n\rfloor$. We have $\alpha(G)<l$ and 
 \beq{GvsT}
 P(K_3,G)-P(K_3,T)\le {2\e n\choose 3}=O(\e^3 n^k),
 \eeq
 as each triangle in $G$ but not in $T$ has to lie inside $U_1\cup U_2$. Let us plug this $G$ into~\req{derivation}. As we have just
observed, the left-hand side of~\req{derivation} is $O(\e^3)$.  Since $G_i$ embeds into
$T'$, we have that $p(G_i,G)\ge \Omega(\e^2)$. (Indeed, if we take any
$f':V(G_i)\to V(G)$ so that $f'(x)$ and $f(x)$ always belong to the same part
of
$\OO T_{l-1}(n)$ and 
$f'(x)\in U_j$ whenever $f(x)=x_j$, then we obtain at least $(1-o(1)) \times (\e
n)^2\times (\frac{n}{l-1})^{k-2}$ different embeddings $f'$.) As $\e$ can be
arbitrarily small, it follows that $G_i$ is sharp by a version of \req{XX}.\epf

Lemma~\ref{lm:phantom} shows that more graphs are necessarily sharp
than those that embed into $\OO T_{l-1}(n)$. Likewise, by unfolding the last inequality in~\req{derivation} and using~\req{GvsT}, we conclude that
$a=O(\e^3 n^N)$. Each of the $t$ summands in 
 \beq{at}
 a=\sum_{i=1}^t\,\sum_{\psi\,:\, G[\psi]\cong \tau_i}
\V x_\psi Q^{\tau_i}\V x_\psi^T
 \eeq 
 is non-negative and is therefore at most $O(\e^3 n^N)$. Thus all terms in
the right-hand side of~\req{at} that can have magnitude $\Omega(\e^2 n^N)$ have to cancel each other. In
particular, for every type $\tau_i$ that embeds into $T'$ but not into
$T$, there are some further zero eigenvectors of $Q^{\tau_i}$ (that are not caught
by the direct application of Lemma~\ref{lm:forced}).

Once we understood ``phantom'' edges, the rounding problem went through without any problems.
The option \verb=phantom_edge= (see the scripts) instructs
\emph{Flagmatic} to take all such extra sharp graphs and zero eigenvectors into
account. 

A similar phenomenon was encountered in the maximum codegree problem for 3-graphs with independent neighbourhoods, see~\cite{falgas+marchant+pikhurko+vaughan}, and a version of
Lemma~\ref{lm:phantom} was crucial for rounding the numerical solution there.

\section{Proving the Stability Property}

Here we prove Theorem~\ref{th:stab}. Our proof is 
similar in 
spirit to the proof of Theorem~2 in \cite{pikhurko:11}. Let $(k,l)$ and $F$ be
as in the theorem. Let $N=N(k,l)$ be the number of vertices
that was used in the flag algebra proof of Section~\ref{lower};
thus $N(3,4)=5$, $N(3,5)=N(4,3)=N(5,3)=6$, $N(3,6)=N(6,3)=7$, and
$N(3,7)=N(7,3)=8$.

Suppose on the contrary that there is
$\e>0$ such that for infinitely many $n\to \infty$
there is a graph $G$ of order $n$ such that
$\alpha(G)<l$ and $p(K_k,G)=c_{k,l}+o(1)$ but
$G$ is $\e {n\choose 2}$-far from a uniform expansion of $F$. 
Let  $V=V(G)$.

Recall that $G_i$ is \emph{sharp} if we have equality in~\req{main}.
Call an admissible graph $G_i$
\emph{singular} if $G_i$ is not contained as an induced subgraph
in any expansion of $F$.   Note that these definitions apply only to
the order-$N$ graphs $G_1,\dots,G_g$.
The following observation is well known (compare it with Lemma~\ref{lm:ForcedSharp}).

\begin{lemma}\label{lm:sharp} Let $i\in [g]$. If $G_i$
is not sharp, then $p(G_i,G)=o(1)$.\end{lemma}

\bpf Note that we have already established that $c_{k,l}'=c_{k,l}$.
Let us run our flag algebra proof on $G$. Similarly to~\req{XX}, 
we obtain that
 $$
 p(K_k,G)-c_{k,l}\ge (p(K_k,G_i)-c_{k,l}-\alpha_i)p(G_i,G)+ O(1/n).
 $$
  Since $G$ is almost extremal, we have that $p(K_k,G)-c_{k,l}=o(1)$. 
The lemma follows from~\req{main}.\epf

\subsection{Cases $(k,l)=(4,3)$ or $(5,3)$}\label{43stab}

Let $l=3$ and $k=4$ or $5$. Here $F$ is the $5$-cycle $C_5$.

The scripts verify that the number of graphs
of order $N$ that occur with positive density in a large expansion of $F$ is the same as the number of sharp graphs (namely, there are 17 graphs in each list). Thus these two lists coincide
by Lemma~\ref{lm:ForcedSharp}. (In other words, each $G_i$ is either sharp or singular.)

By Lemma~\ref{lm:sharp}, we conclude that $p(G_i,G)=o(1)$
for every singular $G_i$. The Induced Removal Lemma of Alon, Fischer,
Krivelevich, and Szegedy~\cite{alon+fischer+krivelevich+szegedy:00} implies that
we can change $o(n^2)$
edges in $G$ and destroy all singular graphs and, additionally, preserve the
property $p(\OO K_3,G)=0$. Since
changing $o(n^2)$ edges affects each $p(H,G)$ by $o(1)$, we can assume that
$G$ itself does not contain any singular induced subgraph. This means the
following.

\Claim1{For any subset $U\subseteq V(G)$ with at most $6$ vertices there is a partition
$U=U_0\cup\dots\cup U_4$ such that $G[U]=\blow{C_5}{U_0,\dots,U_4}$.\ecpf}

By the Induced Removal Lemma we
can additionally assume that either the density of $C_5$ in $G$ is $\Omega(1)$
or $G$ does not have a single induced 5-cycle. In fact, the first alternative necessarily holds:

\Claim{C5}{$p(C_5,G)=\Omega(1)$.}

\bcpf Suppose on the contrary that $G$ does not contain an induced
pentagon.
Take a longest induced path $(u_1,\dots,u_s)$. 
By Claim~\ref{cl:1}, we have $s\le 4$. Also, $s\ge 3$ for otherwise $G$ is
the union of disjoint cliques, of which there can be at most two
because the independence number is at most 2; but then
the $K_k$-density is at least $1/2^{k-1}+o(1)$, contradicting
the extremality
of $G$. Take any vertex $x\in V(G)$. The set $X=\{u_1,\dots,u_s,x\}$
induces some expansion of $C_5$ by Claim~\ref{cl:1}.
Since we do not have an induced pentagon and $s$ is maximal,
$X$ in fact induces an
expansion of the $s$-vertex path $P_s$. Let $\{x,u_i\}$
be the part of this expansion that contains $x$. We assign this vertex 
$x$ into
the $i$-th part, thus obtaining a partition $V(G)=U_1\cup\dots\cup
U_s$.

We have in fact $G=\blow{P_s}{U_1,\dots,U_s}$. Indeed, if we take any two
vertices $x,y$ and apply Claim~\ref{cl:1} to $\{u_1,\dots,u_s,x,y\}$, 
we see that the adjacency relation between $x$ and $y$ in $G$ is exactly as
dictated by the expansion.

Thus we can make $G$ into the union of two disjoint
cliques by removing some edges and without creating $\OO K_3$.
This cannot increase the density of $K_k$ and, as we have just seen,
leads to a contradiction.\ecpf

So suppose that $u_0,\dots,u_4\in V(G)$ span an induced pentagon with
$\{u_i,u_{i+1}\}\in E(G)$
for $i\in \I Z_5$, where $\I Z_5$ denotes the residues modulo $5$. 
Let $U=\{u_0,\dots,u_4\}$.

\Claim{u}{For any $u\in V(G)\setminus U$ there is $j\in \I Z_5$ such that $\{u,u_i\}\in E(G)$
if and only if $i\in \{j-1,j,j+1\}$.}

\bcpf Take the partition $U\cup\{u\}=U_0\cup\dots\cup U_4$ given by
Claim~\ref{cl:1}. For every distinct $i,j\in\I Z_5$, 
the vertices $u_i$ and $u_j$ have different neighbourhoods in
$U\setminus\{u_i,u_j\}$, so
they belong to different parts. Without loss of generality assume that $u_i\in U_i$
for each $i$. If the vertex $u$ belongs to $U_j$, then the neighbours of $u$ are
$u_{j-1},u_j,u_{j+1}$, 
as required.\ecpf 

Claim~\ref{cl:u} gives a partition of $V(G)$ into 5 parts $U_0,\dots,U_4$ where
we classify vertices according to their neighbourhoods in $U$:
 $$
 U_i=\{u_i\}\cup \{u\in V(G)\setminus U\mid \Gamma_G(u)\cap
U=\{u_{i-1},u_i,u_{i+1}\}\}.
 $$

\Claim{Ui}{For every $i\in \I Z_5$ the induced subgraph $G[U_i]$ is complete.} 

\bcpf By symmetry, let $i=0$. Take any distinct $u,v\in U_0$. By the
definition of $U_i$,
we have that $v,u_1,\dots,u_4$ span an induced 5-cycle. Also, $u$ is adjacent
to $u_4$ and $u_1$. By Claim~\ref{cl:u} we conclude that $\{u,v\}\in E(G)$.\ecpf

\Claim{UiUj}{Let $i,j\in \I Z_5$ be distinct and let $v_i\in U_i$ and $v_j\in U_j$
be arbitrary. Then $v_i$ and $v_j$ are adjacent if and only if $i=j\pm 1$.}

\bcpf First, let $i=0$ and $j=1$. The vertex $v_1\in U_1$ is adjacent to
the vertices $u_1$ and $u_2$ but not to $u_3$ of the 5-cycle on $v_0,u_1,\dots,u_4$.
By Claim~\ref{cl:u}, $v_1$ and $v_0$ are adjacent. Next, let $i=0$ and $j=2$. The vertex $v_2\in U_2$ is adjacent
to the vertices $u_1$, $u_2$ and $u_3$ of the 5-cycle on $v_0,u_1,\dots,u_4$.
By Claim~\ref{cl:u}, $v_1$ and $v_0$ are not adjacent. This covers all the cases of Claim~\ref{cl:UiUj} up to a symmetry.\ecpf

Thus we see that $G$ is exactly an expansion of $C_5$ with parts
$U_0,\dots,U_4$, as required.
Choose an arbitrary subsequence of $n$ such that each $|U_i|/n$ approaches
some limit $\alpha_i$.
It remains to show that each $\alpha_i=\frac15$. One approach to
showing this would be to argue that an
explicit degree-$k$ polynomial, that approximates $p(K_k,G)$,
has the unique minimiser $(\frac15,\dots,\frac15)$. This 
approach seems rather messy.

However, there is another way of getting the desired conclusion:
namely, by applying Lemma~\ref{lm:forced}. Let us consider type $\tau_6$ 
which is obtained by labelling the vertices of the 3-edge path
by $3,1,2,4$ as we go along the path. (It is
\texttt{4:121324} in \emph{Flagmatic's} notation.)
There are exactly $8$ non-isomorphic $\tau_6$-flags on $5$ vertices that
we denote by $F_1^{\tau_6},\dots,F_8^{\tau_6}$.
Three of these flags, labelled by \emph{Flagmatic} as
$F_6^{\tau_6},F_7^{\tau_6},
F_8^{\tau_6}$, do not embed into any expansion of $C_5$
when we view them as unlabelled graphs. Thus, by Claim~\ref{cl:1}, 
we have that $p(F_i^{\tau_6},(G,\phi))=0$ for every $\phi$ and $i=6,7,8$.
Every embedding $\psi$ of $\tau_6$ into $G=\blow{C_5}{U_0,\dots,U_4}$ uses
four different parts. The number of embeddings that use a part
of size $o(n)$ is clearly at most $o(n^5)$. So fix an embedding $\psi$
that uses only parts of size $\Omega(n)$.  When we form the vector $\B x_\psi$ as in~\req{x}, we have to
count the number of $\tau_6$-flags on $5$ vertices that we obtain over all
$n-4$ choices of an unlabelled vertex $u\in V(G)\setminus \psi([4])$. Up to symmetry, there are only 5 different choices of $u$ depending on which part $U_i$ contains $u$.
Each $i$ contributes either $|U_i|$ or $|U_i|-1$ to some coordinate of $\B x_\psi$ 
and different $i$'s contribute to different coordinates. Thus, up to a permutation of
coordinates, $\B x_\psi$ is
equal to $(\alpha_1n+o(n),\dots,\alpha_5n+o(n),0,0,0)$.  It follows from a version of Lemma~\ref{lm:forced} that some permutation of
$(\alpha_1,\dots,\alpha_5,0,0,0)$
is a zero eigenvector of $Q^{\tau_6}$. On the other hand, Lemma~\ref{lm:forced} implies that
$(\frac15,\frac15,\frac15,\frac15,\frac15,0,0,0)$ is a forced zero eigenvector of $Q^{\tau_6}$ (that comes from analysing our flag algebra proof on the uniform expansion of $C_5$).
Moreover, the scripts verify 
that the rank of the rational $8\times 8$-matrix $Q^{\tau_6}$ is exactly $7$
(so its null-space has dimension 1).
Since $\alpha_1+\dots+\alpha_5=1$, we conclude that
each $\alpha_i=\frac15$, giving the desired stability property.

\subsection{Cases $k=3$ and $4\le l\le 7$}\label{34stab}

The scripts verify that the number of sharp
graphs and the number of those order-$N$ graphs that embed into
$\OO T_{l-1}(n)$ with one edge added are the same: namely,
10, 20, 33, and  55 graphs when $(l,N)$ is respectively $(4,5)$, $(5,6)$, $(6,7)$, and
$(7,8)$.
Thus these lists coincide by Lemma~\ref{lm:phantom}.
By applying the Induced Removal Lemma, we can assume that
$G$ does not contain any non-sharp $N$-vertex graph. In other words, the
following holds.

\Claim{47}{Every subset $U\subseteq G$ with at most $N$ vertices admits a
partition $U=U_1\cup\dots\cup U_{l-1}$ such that $G[U]$ is equal to
$\blow{\OO K_{l-1}}{U_1,\dots,U_{l-1}}$ with at most one added edge.\ecpf}

Define an equivalence relation $\sim$ on vertices of $G$, 
where $x\sim y$
if and only if $x=y$ or there is a chain of intersecting triangles in $G$ that connects
$x$ to $y$. Each equivalence class is a clique by Claim~\ref{cl:47}
as $N\ge 5$.
Let $U_0$ be the union of equivalence classes of size $1$, that is,
$U_0$ consists of  those vertices that are not contained in a triangle.
Since $G$ does
not contain $\OO K_l$, we have that $|U_0|+1$ is at most the Ramsey number 
$R(3,l)$. Remove $U_0$ from $V(G)$ as this will not affect the stability
property. 

Let $U_1,\dots,U_s$ be the remaining $\sim$-equivalence classes. Each $U_i$
spans a clique and has at least three vertices.

Let us derive a contradiction by assuming that some $U_i$ sends at least two
edges to $V(G)\setminus U_i$, say $\{w,x\}$ and $\{y,z\}$ with $w,y\in U_i$.
Take some $5$-set $X\supseteq \{w,x,y,z\}$ with $|U_i\cap X|=3$. Then $G[X]$
is a subgraph that contains at least one triangle (on $X\cap U_i$) plus at
least two extra
edges. By Claim~\ref{cl:47}, $X$
spans a clique, which contradicts the fact that $x,z\not\in U_i$. 

Thus by removing at most one vertex from each
$U_i$, we can eliminate all edges across the parts. As $U_i$ is still
non-empty, we have that $s< l$ by the $\OO K_l$-freeness of $G$.

A simple optimisation shows
that, in fact, $s=l-1$ and each $U_i$ has $(\frac1{l-1}+o(1))n$ vertices. This
proves the stability property for $f(n,3,l)$ with $4\le l\le 7$.

\subsection{Cases $(k,l)=(6,3)$ or $(7,3)$}

Here $N=7$ if $k=6$ and $N=8$ if $k=7$.
Let $G$ be a $K_3$-free graph
of large order $n$ with $p(\OO K_k,G)=c_{k,l}+o(1)$. Recall that, 
for notational convenience, we prefer to work with
the graph complements in these cases. Also note that an expansion corresponds to a blow-up of a graph when we look at the complements.

The scripts verify that the numbers of the sharp graphs and of those $N$-vertex graphs that appear in a blow-up  of
the Clebsch graph are the same (namely 86 graphs for $(k,N)=(6,7)$ and $232$ graphs
for $(k,N)=(7,8)$).  So these lists coincide by Lemma~\ref{lm:ForcedSharp}. As before,
by applying the Induced Removal Lemma we can additionally assume that
$G$ has the following property:

\Claim{67Sharp}{No singular graph is an induced subgraph of $G$, that is, every 
induced $N$-vertex subgraph of $G$ is a blow-up of  the Clebsch graph $L$.\ecpf}

We need some further definitions before we can proceed with the proof.

Let $X\subseteq V(H)$ be a subset of vertices in some graph $H$. 
Two vertices $x,y\in V(H)$ are \emph{$X$-equivalent}, denoted as $x\sim_X y$,
if $\Gamma_H(x)\cap X=\Gamma_H(y)\cap X$, that is, if they are adjacent to the 
same vertices of $X$. Note that we allow $x$ or $y$ to belong to $X$ and it
is possible that some $x\in X$ and $y\not\in X$ are $X$-equivalent. Clearly,
$\sim_X$ is an equivalence relation. Let $[x]_X=\{y\in V(H)\mid y\sim_X x\}$
denote the equivalence class of $x$. 

Let $C_5'$ be obtained from the 5-cycle on $x_1,\dots,x_5$ by adding an
extra isolated vertex
$x_0$.  Let $\phi$ be a strong homomorphism from $C_5'$ to the
Clebsch graph $L$ that maps the isolated vertex to $00000$ and maps the
remaining
vertices to the cyclic shifts of $00011$. This $\phi$ is injective and its image
is
 \beq{X}
 X=\{00000,00011,01100,10001,00110,11000\}.
 \eeq

\Claim{C5p}{Let $\phi$ and $X$ be as above. Then the following claims hold.
 \begin{enumerate}
 \item For every strong homomorphism $\psi$ of $C_5'$ into
$L$, there is an automorphism $\sigma$ of $L$ such that $\psi=\sigma\circ \phi$.

 \item The $X$-equivalence 
relation is trivial on $V(L)$, that is, $x\sim_X y$ if and only if 
$x=y$. 

\item For every two distinct vertices $x,y\in V(L)$ 
there is $z\in X\setminus\{\phi(x_0)\}$ such that, for $Z=X\setminus\{z\}$,
we have $x\not\sim_Z y$ and the bipartite subgraph of $L$ 
induced by $[x]_Z$ and $[y]_Z$ is either complete or empty.
 \end{enumerate}}

\bcpf Up to an automorphism of $L$, each strong homomorphism $\psi$ from $C_5'$ 
to $L$ is as
follows.
By the vertex-transitivity of $L$, we can assume that $\psi(x_0)=00000$. 
Thus every other vertex of $C_5'$ has to be mapped to a sequence of
weight $2$. (No other vertex can be mapped to $00000$ because $x_0$ 
is the unique isolated vertex of $C_5'$.) 
By permuting indices $1,\dots,5$ (which gives an automorphism of $L$), we can assume that 
$\psi(x_1)=00011$.
Next, up to a permutation of indices $1,2,3$, 
we can assume that $\psi(x_2)=01100$ and $\psi(x_5)=11000$. Up to
a transposition of  $4$ and $5$, we can also assume that $\psi(x_3)=10001$. Now,
$\psi(x_4)=00110$ is uniquely determined.
Thus $\psi=\phi$ up an automorphism of $L$, proving the first part of the
claim. 

Every 5-sequence of weight 0, 4 and 2 sends respectively
0, 3, and 1--2 edges to $X$, so $X$ distinguishes
vertices of different weight. An easy case analysis
for each possible weight shows the second part of the claim. For example,
$00011$ is identified among all weight-2 sequences already by the set
$\{01100,11000\}\subseteq X$.

In order to establish the third part, we use the fact that any 
cyclic permutation or reversal of the indices preserves $X$. Up to
these symmetries, there are 12 different unordered pairs
$x,y$ to check. The following table lists a vertex
$z$ that establishes the claim and
the $Z$-equivalence classes of $x$ and $y$,
where $Z=X\setminus\{z\}$:
 $$
 \begin{array}{c|c|c|c|c}
 x & y & z & [x]_Z & [y]_Z\\ \hline
 00000 & 00011 & 10001 & \{00000,01010\} & \{00011\}\\
 00000 & 00101 & 00011 & \{00000,10100\} & \{00101\}\\
 00000 & 01111 & 00011 & \{00000,10100\} & \{01111\}\\
 00011 & 01100 & 00110 & \{00011\} & \{ 10100\}\\
 00011 & 00110 & 00011 & \{00011\} & \{00110\}\\
 00011 & 00101 & 00011 & \{00011\} & \{00101\}\\
 00011 & 01010 & 00110 & \{00011\} & \{01010\}\\
 00011 & 10100 & 10001 & \{00011\} & \{01100,10100\}\\
 00101 & 01010 & 00110 & \{00101\} & \{01010\}\\
 00101 & 01001 & 00110 & \{00101\} & \{00000,01001\}\\
 01111 & 10111 & 00011 & \{ 01111\} & \{10111\}\\
 01111 & 11011 & 00011 & \{01111\} & \{11011\}
 \end{array}
$$
Alternatively, the included \emph{Mathematica} 
notebook \texttt{Clebsch.nb} available from the ancillary folder
of~\cite{pikhurko+vaughan:Fkl:arxiv}  
verifies the existence of $z$ by
the brute-force enumeration of all cases. This proves Part 3 of the claim.\ecpf

\Claim{C5P}{$\pp(C_5',G)=\Omega(n^6)$.}
\bcpf  Suppose on the contrary that $p(C_5',G)=o(1)$. By the Induced Removal Lemma,
we can additionally assume that $\pp(C_5',G)=0$. We let \emph{Flagmatic} prove some lower bound
on the density of $\OO K_k$ given that both $K_3$ and $C_5'$ are forbidden.
The obtained bound (with the certificates \texttt{63a.js}
and \texttt{73a.js}) is strictly larger
than $c_{k,3}$. This contradicts $p(\OO K_k,G)=c_{k,3}+o(1)$ for all large $n$, 
proving the claim.\ecpf

Fix one embedding $\psi$ of $C_5'$ into $G$. Let us view $C_5'$ as the
subgraph of $L$ induced by $X\subseteq V(L)$, where $X=V(C_5')$ is defined
by~\req{X}. Thus $\psi:X\to V(G)$. Let $Y=\psi(X)$. 

\Claim{16}{For every $y\in V(G)$ there is the (unique) vertex $x\in V(L)$ whose
adjacencies to $X$ match those of $y$ to $Y$, that is, $\psi(\Gamma_L(x)\cap
X)=\Gamma_G(y)\cap Y$.}

\bcpf The subgraph $H=G[Y\cup\{y\}]$, that has
at most $7\le N$ vertices, admits an embedding 
into a blow-up of the Clebsch graph by Claim~\ref{cl:67Sharp}. This implies
that there is a strong homomorphism $\xi$ from $H$ into $L$. By Part 1 of
Claim~\ref{cl:C5p}, we can assume that the composition $\xi\circ \psi$ 
is the identity map $\mathrm{Id}_X: X\to X$. Now, $x=\xi(y)$ satisfies the claim. The uniqueness of $x$ follows from Part 2
of Claim~\ref{cl:C5p}.\ecpf

Thus each $y\in V(G)$ falls into one of at most sixteen $Y$-equivalence classes that are naturally labelled as $U_x$ for $x\in V(L)$, where $x=x(y)$ is given by Claim~\ref{cl:16}. In particular, for each $x\in X$, 
the part containing $\psi(x)$ is labelled by $U_x$.

\Claim{Between}{For every adjacent $x,y\in V(L)$, the induced bipartite subgraph $G[U_x,U_y]$ is
complete. For non-adjacent $x,y\in V(L)$ the induced bipartite
subgraph $G[U_x,U_y]$ is empty. (In particular, each part $U_x$
forms an independent
set.)}

\bcpf Let $x,y\in V(L)$ be adjacent. Let $x'\in U_x$ and $y'\in U_y$ be arbitrary. 

Pick $z\in X$ given by Part~3 of Claim~\ref{cl:C5p} and let $Z=X\setminus\{z\}$.
The induced subgraph $H=G[\psi(Z)\cup\{x',y'\}]$ has at most $7\le N$ vertices.
By Claim~\ref{cl:67Sharp}, $H$ admits a strong homomorphism $\xi$ to
$L$. 
By Part 1 of Claim~\ref{cl:C5p}, 
we can assume that $\xi\circ \psi$ is the identity on $Z$. Then $\xi(x')\in [x]_Z$ and $\xi(y')\in [y]_Z$. However, the bipartite
subgraph induced by $[x]_Z$ and $[y]_Z$ in $L$ is complete by the choice
of $z$ (since $\{x,y\}\in E(L)$). Thus $x'$ and $y'$ are adjacent.
The second part of the claim follows in a similar manner.\ecpf

Thus we know that $G$ is a blow-up of $L$ with parts $U_{00000},\dots,U_{11110}$. It remains to argue that
each part $U_x$ has $(\frac1{16}+o(1))n$ vertices. 

Let $k=7$. We proceed very similarly
as we did at the end of Section~\ref{43stab} so we are rather brief.
We consider the type $\tau_{37}$, which is a labelling of $C_5'$.
It is \texttt{6:1213243545} in \emph{Flagmatic's} notation. 
There are 22 $\tau_{37}$-flags on 7 vertices. By Claim~\ref{cl:C5P},
there are  $\Omega(n^6)$ embeddings $\psi$ of $\tau_{37}$ into $G$.
By Parts 1--2 of
Claim~\ref{cl:C5p}, each obtained vector $\B x_{\psi}$ consists of 
sixteen entries $|U_x|+O(1)$, one for each $x\in V(L)$,
and six zeros.  On the other hand, the script \texttt{73.sage} verifies that
the $22\times 22$-matrix
$Q^{\tau_{37}}$ from our flag algebra proof has rank $21$.
Moreover,  by Lemma~\ref{lm:forced}, the matrix $Q^{\tau_{37}}$ has one forced zero eigenvector 
consisting of $16$ entries equal to $1/16$ and
six entries equal to $0$.  It follows in the same way as in  Section~\ref{43stab} that
each $U_x$ has size $(\frac1{16}+o(1))n$.

Let $k=6$. We consider the
type $\tau_{11}$ that consists of the $3$-edge path plus an 
isolated vertex (it is \texttt{5:121324} in \emph{Flagmatic's} notation). 
Since
$C_5'$ contains $\tau_{11}$ as a subgraph, Claim~\ref{cl:C5P}
implies that there are $\Omega(n^5)$ embeddings $\xi$ of $\tau_{11}$ into $G$. Fix an embedding $\xi$ such that its image avoids all parts $U_x$ of size $o(n)$. (A typical $\xi$ has
this property.) Similarly to Part 1 of Claim~\ref{cl:C5p}, we can relabel the parts $U_x$ so that 
the image $Y$ of $\xi$ has exactly one vertex in each of the parts $U_{00000},U_{00011},U_{01100},U_{10001},U_{00110}$. The $Y$-equivalence
relation on $G$ makes each part $U_x$ into a separate equivalence class except for the following three $Y$-equivalence classes:
 \beq{exception}
 U_{00000}\cup U_{00101},\quad U_{00011}\cup U_{10010},\quad U_{00110}\cup 
U_{01010}.
 \eeq
 On the other hand, the $16\times 16$-matrix $Q^{\tau_{11}}$ of our solution has rank $15$. Moreover, it
has one forced zero eigenvector that has $10$ entries
equal to $1/16$ and $3$ entries equal to $2/16$ by Lemma~\ref{lm:forced}. (This follows from~\req{exception} when applied to the uniform blow-up of $L$.) 
This implies that each of the 10 parts that do no appear in~\req{exception}
has size $(\frac1{16}+o(1))n$ while each of the three sets in~\req{exception} has $(\frac2{16}+o(1))n$ vertices. 

The graph $G$ has other copies  of $\tau_{11}$, e.g.\ via $$
 U_{10100},U_{01111},U_{11000},U_{10111},U_{11101}.
 $$
 The adjacency pattern to these $(\frac{n}{16}+o(n))^5$ copies $\tau_{11}$ uniquely
identifies parts $U_{00000}$,
$U_{00101}$, $U_{00011}$ and $U_{01010}$. As before, we conclude that that each of these
parts has size $(\frac1{16}+o(1))n$. This is enough to determine the
sizes of those parts that appear in~\req{exception}. Thus $G$ is $o(n^2)$-close to a uniform blow-up of
$L$. The stability property has been established.

\brm By running
everything with $N=8$ (see the script \texttt{63.sage} and the certificate \texttt{63b.sage}), it is possible to shorten the ``human'' part of the proof
of Theorem~\ref{th:stab} for $(k,l)=(6,3)$. (Namely, Part~3 of
Claim~\ref{cl:C5p} and the argument around \req{exception} become redundant.) However,
we believe that the ability to solve this case within the universe of 7-vertex
graphs justifies the extra work, as the ideas introduced for this task may be
useful for other problems.\medskip

\section{Exact Result}

First,  we present a rather general Theorem~\ref{th:1} and then verify in
Section~\ref{>exact} that it implies Theorem~\ref{th:exact}. 
Theorem~\ref{th:1} could in principle be strengthened in various ways
but we state only the current version as it suffices for all
the cases that we need.

\subsection{A General Result}

We need to give some definitions first, given an arbitrary pair $(k,l)$ and any 
admissible graph
$F$ with vertex set $[m]$.

We say that $F$ is a \emph{stability graph}
for $(k,l)$ if
for every $\e>0$ there are $n_0$ and $\delta>0$ such that the following holds.
Let $G$ be an arbitrary graph such that $n=v(G)\ge n_0$,
$\alpha(G)<l$, and $p(K_k,G)\le c_{k,l}+\delta$. 
Then there is a partition $V(G)=V_1\cup\dots\cup V_m$ such that the part
sizes differ at most by 1 and
 $$
 |E(\blow{F}{V_1,\dots,V_m})\bigtriangleup E(G)|\le \e {n\choose 2}.
 $$
 In other words, $F$ is a stability graph for $(k,l)$ if every
large almost extremal graph for the $f(n,k,l)$-problem is $o(n^2)$-close in
the edit distance to a uniform expansion of $F$. Clearly, this property is
preserved if we replace $F$ by an isomorphic graph or
by $\blow{F}{U_1,\dots,U_m}$ with $|U_1|=\dots=|U_m|>0$. 

We give some further definitions related to the graph $F$, which will be illustrated in
the next
paragraph. Let us call a set of vertices
$X\subseteq[m]$
\emph{legal} if $F-X$ does not contain $\OO K_{l-1}$. Let the 
\emph{gradient} $\grad(X)$ of $X$ be the probability, when we pick $k-1$
independent and uniformly distributed vertices $x_1,\dots,x_{k-1}\in [m]$,
that all belong to $X$ and for every $i,j\in [k-1]$ the vertices
$x_i$ and
$x_j$ are adjacent or equal. 
Let us call a stability graph $F$ \emph{strict} 
if $\grad(X)>c_{k,l}$ for every legal $X$
for which there is no $i\in [m]$ with $X=\hat\Gamma_F(i)$. Recall that 
 $$
 \hat\Gamma_F(i)=\{i\}\cup \{j\in V(F)\mid \{i,j\}\in E(F)\}
 $$
 is the closed neighbourhood of $i$.

The above definitions are
motivated by
the addition of
a new vertex $x$ to $F'=\blow{F}{V_1,\dots,V_m}$ with $|V_1|=\dots=|V_m|=n/m$
so that $x$ is adjacent to precisely $\cup_{i\in X} V_i$. The new graph
is still $\OO K_l$-free if and only if $X$ is legal. Also, the number
of $k$-cliques that contain $x$ is $\grad(X){n\choose k-1}+O(n^{k-2})$.
If $X=\hat\Gamma_F(i)$, then adding $x$ is the same as enlarging the
part $V_i$ by one vertex and, if $F$ is a stability
graph, then the number of $k$-cliques increases by 
$(c_{k,l}+o(1)){n\choose k-1}$, see Claim~\ref{cl:GKkReg} below. 
Thus $F$ is
strict if the number of the new $k$-cliques is by $\Omega(n^{k-1})$ larger for
every other legal $X$.

\begin{theorem}\label{th:1} Let a pair $(k,l)$ admit a stability graph $F$
which is strict. Then there is $n_0$ such that every graph $G$ with
$n=v(G)\ge n_0$, $\alpha(G)<l$, and $P(K_k,G)=f(n,k,l)$ contains an expansion
of $F$ as a spanning subgraph.
\end{theorem}

\bpf Let $V(F)=[m]$. Choose positive constants
 \beq{constants}
 \e_2\gg \e_1\gg \e_0\gg
1/n_0>0,
 \eeq 
 each being sufficiently small, depending on the previous ones.
We show that $n_0$ satisfies the conclusion of the theorem.

Since there are finitely many different subsets
$X\subseteq[m]$, we can assume that 
 \beq{grad}
 \grad(X)\ge c_{k,l}+2km\e_2
 \eeq
 for every legal $X$ that is not the closed neighbourhood of some vertex.
Also,
we may assume that for every $n\ge n_0$ we have 
 \beq{rate}
 f(n,k,l)\ge (c_{k,l}-\e_0){n\choose k},
 \eeq

Let $G$ be an arbitrary $f(n,k,l)$-extremal graph with $n\ge n_0$
vertices. Let $V=V(G)$. 
Since $f(n,k,l)=(c_{k,l}+o(1)){n\choose k}$
by~\req{ckl} and $F$ is a stability graph,
we have that
 \beq{F0'}
 |E(G)\bigtriangleup E(F')|\le \e_0{n\choose 2}
 \eeq 
 for some
uniform
expansion $F'=\blow{F}{V_1,\dots,V_m}$ on $V$.

We are going to modify the partition $V=V_1\cup\dots\cup V_m$.
Given a current partition, let $B=E(F')\setminus E(G)$ and $S=E(G)\setminus
E(F')$. We call the pairs in $B$ \emph{bad} and those in $S$ 
\emph{superfluous}. 

Iteratively repeat the following operation as long as possible (updating
$V_1,\dots,V_m$, $F'$, $B$ and $S$ as we proceed):
if we can move some vertex $x$ of $F'$ to another part and decrease
the number of bad pairs by least $\e_1 n$, then we perform this move.

Since we had initially at most $\e_0{n\choose 2}$ bad pairs, we perform at most
$\e_0{n\choose 2}/\e_1 n <\e_1 n/4$ moves. Let $V_1,\dots,V_m,F',B,S$
refer to
the final configuration.
What we have achieved is that for every vertex $x\in V_i$ and every $j\in [m]$
 \beq{move}
 |\Gamma_{\OO G}(x)\cap \cup_{h\in\hat\Gamma_F(j)}V_h| > |\Gamma_{\OO G}(x)\cap \cup_{h\in\hat\Gamma_F(i)}V_h|-\e_1n.
 \eeq 
 Also, the current expansion $F'$ is not far from being uniform:
 \beq{Vi}
 \left|\,|V_i|-\frac nm\,\right|\le \e_1n, \quad \mbox{for all
$i\in [m]$}.
 \eeq
 In addition, we have
 \beq{e1}
 \left| E(G)\bigtriangleup E(F')\right| \le \e_0{n\choose 2}+ \frac{\e_1 n}4 \,
n< \e_1{n\choose 2}.
 \eeq

\Claim{OKl}{The removal of any edge $\{x,y\}$ from $F'$ creates $\OO K_l$.}

\bcpf First, suppose that $x$ and $y$ belong the same part $V_i$.
Partition $V_i=X\cup Y$ into two almost equal parts so that
$x\in X$ and $y\in Y$. Let
$F''$ be obtained from
$F'$ by removing all edges between $X$ and $Y$. 
By~\req{Vi} and~\req{e1} we have rather roughly that
 \begin{eqnarray*}
 \pp(K_k,F'')&\le& \pp(K_k,F')- \frac1{2^{k}} {|V_i|\choose k}\\
  &\le& \pp(K_k,G)+\e_1
{n\choose 2} {n-2\choose k-2} - \frac{(n/m)^k}{2^{k+1} k!}\ <\ \pp(K_k,G).
 \end{eqnarray*}
  By the extremality of $G$, we conclude that $F''$ contains
an independent set $I$ of size $l$. Clearly,
$I$ has exactly one vertex in each $X$ and $Y$. Since any permutation
of the vertices of $X$ (and of $Y$) is an automorphism of $F''$, we can
assume that $x,y\in I$,
giving the required. 

If $x,y$ come from different parts
$V_i$ and $V_j$, then a similar argument works where we
remove all edges of $F'$ between $V_i$ and $V_j$.\ecpf

\Claim{dS}{For every bad pair $\{x_1,x_2\}\in B$ we have $d_S(x_1)+d_S(x_2)\ge
n/(3m^{l-2})$.}

\bcpf Let $x_1\in V_{i_1}$ and $x_2\in V_{i_2}$. By Claim~\ref{cl:OKl},
$F'-\{x_1,x_2\}$ has $\OO K_l$ as a subgraph. 
This means that we can find
distinct
$i_3,\dots,i_k\in [m]\setminus\{i_1,i_2\}$ such that no pair of vertices
$i_1,\dots,i_l$, except $\{i_1,i_2\}$, is adjacent in $F$. 

For every choice of $\B x=(x_3,\dots,x_l)$ such that $x_j\in V_{i_j}$, at
least one pair $\{x_j,x_h\}$ with $1\le j< h\le l$ 
is superfluous (for otherwise
we get an
independent set of size $l$ in $G$). It is impossible that both $j$ and $h$
are at least $3$ for at least half
of the choices of $\B x$: otherwise, as each superfluous pair is overcounted
at most $n^{l-4}$ times, we would have that
 $$
 |S|\ge \frac12 \left(\left(\frac1m-\e_1\right)n\right)^{l-2}  \frac{1}{n^{l-4}}>
\e_1{n\choose 2},
 $$
 which contradicts \req{e1}. Thus, for at least half of the choices of $\B x$
there is a superfluous pair intersecting $\{x_1,x_2\}$. Since each such 
pair is over-counted at most $n^{l-3}$ times, we obtain that
 $$
 d_S(x_1)+d_S(x_2)\ge \frac 12 \left(\left(\frac1m-\e_1\right)n\right)^{l-2}
\times \frac1{n^{l-3}},
 $$
 which implies the claim provided that $\e_1=\e_1(m,l)$ is sufficiently
small.\ecpf

Let $K_k^1$ be the flag obtained from $K_k$ by labelling one vertex. 
Thus $\pp(K_k^1,(H,x))$ is the number of $k$-cliques in a graph 
$H$ that contain $x\in V(H)$.

\Claim{GKkReg}{For any two vertices $x,y\in V$, we have 
 $$
 \left|\pp(K_k^1,(G,x))-\pp(K_k^1,(G,y))\right|\le {n-2\choose k-2}.
 $$}
 \bcpf If we delete $x$ but add a clone $y'$ of $y$ (putting
an edge between $y$ and $y'$), then we do not create a copy of
$\OO K_l$ while the number of $k$-cliques changes by at most
$\pp(K_k^1,(G,y))-\pp(K_k^1,(G,x))+{n-2\choose k-2}$. Since $G$ is extremal, this
has to be non-negative. By swapping the roles of $x$ and $y$, we derive the
claim.\ecpf

 Claim~\ref{cl:GKkReg} and the extremality of $G$ imply that for every $x\in V(G)$ we have
 \beq{rate1}
 \pp(K_k^1,(G,x))\le \frac{k\, f(n,k,l)}n+{n-2\choose k-2},
 \eeq 
 for otherwise $\pp(K_k,G)=\frac1k \sum_{y\in V(G)} \pp(K_k^1,(G,y))>
\frac nk(\pp(K_k^1,(G,x))-{n-2\choose k-2})$ is too large. 

Suppose that $B$ is not empty for otherwise we are done: 
$G$ contains $F'$ as a spanning subgraph.

By Claim~\ref{cl:dS}, there is a vertex $x$ whose $S$-degree is at least
$n/6m^{l-2}$. Define
 $$
 X=\{i\in [m]\mid |V_i\setminus \Gamma_G(x)|\le \e_2 n\}.
 $$

\Claim{legal}{$X$ is legal.}

\bcpf Suppose that this is false. Then there are
distinct $i_1,\dots,i_{l-1}\in
[m]\setminus X$ that span $\OO K_{l-1}$ in $F$. Let $x_l=x$. For every choice
of $(x_1,\dots,x_{l-1})$ with $x_j\in \Gamma_{\OO G}(x)\cap V_{i_j}$, the
$(l-1)$-set
$\{x_1,\dots,x_{l-1}\}$ has to span at least one edge in $G$ (otherwise together
with $x$ it induces $\OO K_l$). This edge is necessarily in $S$. On the other
hand,
any pair in $S$ is over-counted at most $n^{l-3}$ times. Thus
$|S|\ge (\e_2n)^{l-1}/n^{l-3}$, contradicting~\req{e1}.\ecpf

\Claim{X=Gamma}{There is $i\in [m]$ such that $X=\hat\Gamma_F(i)$.}

\bcpf Suppose that the claim is false. As $F$ is strict, we have that \req{grad}
holds. Let $F''$ be obtained from $F'$ by changing edges at $x$ so that the new neighbourhood of $x$ is exactly $Y=(\cup_{j\in X} V_j)\setminus \{x\}$. The number
of $K_k$-subgraphs in $F''$ via $x$ is 
 \beq{F''}
 \pp(K_k^1,(F'',x))\ge (c_{k,l}+2km\e_2){n-1\choose k-1}-\e_1 m n {n-2\choose k-2}+O(1/n).
 \eeq
  (Here, the middle term corresponds to the fact that, by~\req{Vi}, we can make $F'$ into a uniform expansion by moving at most $\e_1 mn$ vertices between parts.) On the other hand, $G$ and $F'$
differ in at most $\e_1{n\choose 2}$ edges by~\req{e1} while at most $\e_2m n$ edges
between $x$ and $Y$ can be missing in $G$ by the definition of $X$. Thus, rather roughly,
 $$
 \pp(K_k^1,(G,x)) \ge \pp(K_k^1,(F'',x))
 - \e_1{n\choose2} {n-3\choose k-3} - \e_2mn{n-2\choose k-2}.
 $$
 However, this 
inequality contradicts \req{rate}, \req{rate1} and \req{F''} by our
choice of the constants in~\req{constants}.\ecpf

Fix $i$ that is returned by Claim~\ref{cl:X=Gamma}.

\Claim{smallDeltaB}{$d_B(x)< 2\e_1 n$.}

\bcpf Suppose on the contrary that $d_B(x)\ge 2\e_1 n$. 

Consider moving
$x$ to $V_i$. (The following statements are also true if $x$ is already 
in $V_i$.) 
By~\req{move}, the new number of bad
pairs at $x$ would be at least $d_B(x)-\e_1 n> \e_2m n$ and each one
would connect $x$
to $\cup_{h\in \hat\Gamma_F(i)} V_h$.
 
Hence, in the graph $G$,  $x$ has more than $\e_2n$ non-neighbours in some $V_h$ with $h\in
\hat\Gamma_F(i)$,
meaning that $X\not=\hat\Gamma_F(i)$ and contradicting Claim~\ref{cl:X=Gamma}.\ecpf

Let $x\in V_j$ (where possibly $j=i$). Fix $y\in V_j$ that has at most the average number
of superfluous edges over the vertices of $V_j$. We have
 $$
 d_S(y)\le \frac{|E(G)\bigtriangleup E(F')|}{|V_j|}\le \frac{\e_1{n\choose 2}}{(1/m-\e_1)n}\le \e_1mn.
 $$
 This and Claim~\ref{cl:smallDeltaB} imply that 
 $$
 |\Gamma_G(y)\setminus \Gamma_G(x)|\le d_S(y)+d_B(x)\le\e_1(m+2)n.
 $$ 
 On the other hand, $x$ sends at least $d_S(x)/m\ge n/6m^{l-1}$
superfluous edges to some part $V_h$. By~\req{e1}, all but
at most $\e_1{n\choose 2}$ pairs of $V_h$ are edges of $G$. Thus the superfluous
edges at $x$ create at least 
 $$
 {n/6m^{l-1}\choose k-1}- \e_1{n\choose 2} {|V_h|-2\choose k-3} >  
(2m+5)\e_1n
{n-2\choose k-2}
 $$
 copies of $K_k$ through $x$. We conclude that
 $$
 \pp(K_k^1,(G,x))-\pp(K_k^1,(G,y))> (2m+5)\e_1n
{n-2\choose k-2} - 2\e_1(m+2)n{n-2\choose k-2}={n-2\choose k-2},
 $$
 contradicting Claim~\ref{cl:GKkReg}.
This final contradiction to $B\not=\emptyset$ proves
Theorem~\ref{th:1}.\qed

\subsection{Verifying Theorem~\ref{th:exact}}\label{>exact}

Theorems~\ref{th:stab} and \ref{th:1} imply Theorem~\ref{th:exact} provided
we can verify that the appropriately defined $F$ is strict. The cases $F=\OO K_{l-1}$
or $C_5$ are straightforward to verify. Namely, every legal set $X$ that is not a closed neighbourhood of a vertex has at least $2$ vertices for $\OO K_{l-1}$ and at least 4 vertices
for $C_5$; any such $X$ contains some closed neighbourhood as a proper subset
and has a strictly larger gradient.

Let $(k,l)=(6,3)$ or 
$(7,3)$. Let us check that $\OO L$ satisfies Theorem~\ref{th:1}.
We already know by Theorem~\ref{th:stab}  that 
$\OO L$ is a stability graph for $(k,l)$.
 Let $X\subseteq
V(L)$ be any legal set, meaning that $Y=V(L)\setminus X$ spans no edge
in $L$. By the 
vertex-transitivity of $L$, we can assume that $00000\in Y$. Thus all
other sequences in $Y$ have weight 2 and, furthermore, no two such
sequences can have 1s in disjoint positions. If $|Y|=5$, then up to a symmetry
the only possibility is $Y=\{00000,00011,00101,01001,10001\}$ but then $X$ is
precisely the closed neighbourhood of $11110$ in $\OO L$. If $|Y|=4$ and $X$
does
not contain a closed neighbourhood, then, up to an automorphism of $L$,
we have $Y=\{00000,00011,00101,00110\}$.
The script \texttt{Clebsch.nb} shows that, if $k=6$, then
$\grad(X)=1437/2^{16}>c_{6,3}$ and if $k=7$, then
$\grad(X)=14503/2^{21}>c_{7,3}$. 
Every other $Y$ is a subset of one
of the sets that we have already considered and the
gradient of $X=V(L)\setminus Y$ is strictly larger than what we had before. Thus 
$\barL$ is strict. This finishes the remaining cases of Theorem~\ref{th:exact}.

\section{Concluding Remarks}

Let us call a graph $G$ \emph{extremal
$(s,t)$-Ramsey} if $G$ has neither $K_s$ nor $\OO K_t$
as an induced subgraph while the order of $G$ is $R(s,t)-1$,
that is, maximum possible. Das et al~\cite{das+huang+ma+naves+sudakov:12:arxiv}
asked if for every $(k,l)$ and large $n$, the value of $f(n,k,l)$ is attained by 
an expansion of some extremal Ramsey graph. The cases $(k,l)=(6,3)$ and $(7,3)$ that we solved here
show that the answer is in the negative. Interestingly,
$\barL$ is nonetheless related to Ramsey numbers, but to 3-colour ones:
Kalbfleisch and Stanton~\cite{kalbfleisch+stanton:68} showed that
there are two different 3-edge-colourings of $K_{16}$ without a monochromatic
triangle but each colour class (in either colouring) is isomorphic to the
Clebsch graph (and thus the union of any two colour classes is isomorphic to $\barL$).

Das et al~\cite[Section 6]{das+huang+ma+naves+sudakov:12:arxiv} mention that they ran
the SDP-solver for the cases $(k,l)=(5,3)$, $(3,5)$ and $(3,6)$ and the obtained floating-point
bound suggested that $c_{5,3}=31/625$, $c_{3,5}=1/16$, and $c_{3,6}=1/25$ with extremal
configurations being an expansion of respectively $C_5$, $\OO K_4$ and $\OO K_5$. Since
their paper was already quite long they did not try to convert it into a rigorous proof.
The current paper makes these statements rigorous.

It would be
interesting to identify further pairs $(k,l)$ amenable to this approach. One promising case is $f(n,4,4)$, where we make the following conjecture.

\begin{conjecture}\label{cj:44}
 \beq{cj:44}
 c_{4,4}=\frac{-11+14\times 2^{1/3}}{192}.
 \eeq
\end{conjecture}

The upper bound in~\req{cj:44}
comes from taking expansions of the (unique) 
$(3,4)$-Ramsey graph $F$ with 8 vertices with 10 edges. More specifically, 
let $F$ be obtained from the $8$-cycle on $1,\dots,8$ by adding the two "diameters" $\{1,5\}$ and $\{2,6\}$ as edges. Take an expansion $F'=\blow{F}{U_1,\dots,U_8}$ with parts $U_1$, $U_2$, $U_5$, and $U_6$
(those corresponding to degree-3 vertices of $F$) having size $(\alpha +o(1))n$ and the other four parts having size $(\frac14-\alpha+o(1))n$, where $\alpha=\frac1{12}+2^{1/3}-2^{2/3}$. Routine calculations show that the density of $K_4$ approaches the right-hand side of~\req{cj:44} as
$n\to\infty$. On the other hand, \emph{Flagmatic} suggests that this construction is asymptotically optimal and, perhaps, a flag algebra proof exists within the 8-vertex universe (i.e.\ taking $N=8$). Unfortunately, we have not been able to round the floating point solution.

\section*{Acknowledgements}

The authors are grateful to the anonymous referee for the careful reading and numerous helpful remarks.

\end{document}